\documentclass[12pt]{amsart}
\usepackage{amsfonts}
\newcommand{\sgn}{\mathop{\rm sgn}}
\newtheorem{thm}{Theorem}

\newcommand{\R}{{\mathbb R}}
\newcommand{\Z}{{\mathbb Z}}

\newcommand{\D}{\Delta}
\newcommand{\Om}{\check{\omega}}
\renewcommand{\u}{\vec u}
\newcommand{\sumkinf}{\sum_{k=0}^\infty}
\newcommand{\sumtinf}{\sum_{t=0}^\infty}
\newcommand{\sumweyl}{\sum_{w\in W}}

\newcommand{\prodin}{\prod_{i=1}^n}
\newcommand{\gel}{g_{\eta\lambda}}
\newcommand{\gpel}{g'_{\eta\lambda}}
\newcommand{\belk}{b_{\eta\lambda,k}}
\newcommand{\bpelk}{b'_{\eta\lambda,k}}
\newcommand{\cgamk}{c_{\gamma, k}}
\newcommand{\cpgamk}{c'_{\gamma, k}}
\def\1/2{\frac{1}{2}} 
\title[Random Walk in an Alcove of an Affine Weyl Group]{Random Walk in
an Alcove of an Affine Weyl Group, and Non-Colliding Random Walks on an
Interval}
\author{David J. Grabiner}\thanks{Supported by an NSF Postdoctoral
Fellowship.} 
\begin{document}
\maketitle
{\bf Abstract.}  We use a reflection argument, introduced by Gessel and
Zeilberger, to count the number of $k$-step walks between two points
which stay within a chamber of a Weyl group.  We apply this technique to
walks in the alcoves of the classical affine Weyl groups.  In all cases,
we get determinant formulas for the number of $k$-step walks.  One
important example is the region $m>x_1>x_2>\cdots>x_n>0$, which is a
rescaled alcove of the affine Weyl group $\tilde C_n$.  If each
coordinate is considered to be an independent particle, this models $n$
non-colliding random walks on the interval $(0,m)$.  Another case models
$n$ non-colliding random walks on a circle.

\section{Introduction}

The {\em ballot problem}, a classical problem in random walks, asks how
many ways there are to walk from the origin to a point
$(\lambda_1,\ldots,\lambda_n)$, taking $k$ unit-length steps in the
positive coordinate directions while staying in the region $x_1\ge
x_2\ge\cdots\ge x_n$.  The solution is known, and leads to the
hook-length formula for Young tableaux; a combinatorial proof, using a
reflection argument, is given in~\cite{WM,Z}.

The same reflection argument has also been applied to the case of $n$
independent diffusions, or $n$ discrete processes which cannot pass each
other without first colliding.  Using this method, Karlin and
McGregor~\cite{KM,KM2} give a determinant formula for the probability or
measure for the $n$ particles, starting at known positions, not to have
collided up to time $t$ and to be in given positions.  Hobson and
Werner~\cite{HW} generalize this argument to $n$ independent Brownian
motions in an interval or on a circle.

Gessel and Zeilberger~\cite{GZ}, and independently Biane~\cite{biane92},
consider a more general question, for which some of the same techniques
apply.  For certain ``reflectable'' walk-types, we can count the number
of $k$-step walks between two points of a lattice, staying within a Weyl
chamber.  The argument generalizes naturally to ``reflectable''
diffusions~\cite{brownian}; we can compute the density function for the
diffusion started at a point $\eta$ to stay within the Weyl chamber up
to time $t$ and be at a point $\lambda$.

Grabiner and Magyar~\cite{decomp} classify all reflectable random walks
for finite Weyl groups, and compute determinant formulas for many
important cases, including walks with steps $\pm e_i$ or with steps
$\pm\1/2 e_1\pm\cdots\pm\1/2 e_n$ in all of the classical Weyl groups.

We prove analogous results for the alcoves of affine Weyl groups.  In
contrast to the chambers of classical Weyl groups, these are bounded
regions, such as $m>x_1>\cdots>x_n>0$.  The reflectable random walks for
the affine Weyl groups are the same as for the corresponding classical
Weyl groups.  We use these reflection arguments to find determinant
formulas for the number of walks of length $k$ which stay within the
alcoves of the classical affine Weyl groups.  We then simplify the
determinant of infinite sums to get a determinant of finite sums of
sines and cosines or of exponentials of cosines, depending on the random
walk.

Many results are known in the $\tilde A_{n-1}$ cases, in which the
region in $\R^2$ is $x_1>x_2>\cdots>x_n>x_1-m$.  The $\tilde A_1$ case
(or equivalently $\tilde B_1$, region $m>x>0$) is a single random walk
on an interval.  This is the classical gambler's ruin problem; gamblers
with initial stakes of $\eta$ and $N-\eta$ chips bet one chip at a time
until one is broke.  The probability that the gambler who started 
with $\eta$ will first go broke after $k$ bets is
\begin{equation}\label{onedim}
\frac{1}{N}\sum_{r=1}^{N-1}
\cos^{k-1}(\pi r/N)\sin(\pi r/N)\sin(\pi\eta r/N).
\end{equation}
This formula goes back to Lagrange~\cite[p.\ 353]{feller}.  Similar
calculations show that the probability that the gambler who started with
$\eta$ will have $\lambda$ left after $k$ bets, with neither player
going broke, is 
\begin{equation}\label{onedimgen}
\frac{2}{N}\sum_{r=1}^{N-1}
\cos^{k}(\pi r/N)\sin(\pi\lambda r/N)\sin(\pi\eta r/N).
\end{equation}
The reflection principle was applied to this problem by
Grossman~\cite{grossman}.  A $q$-analogue of this formula was computed
by Krattenthaler and Mohanty~\cite{KrM}.

The $n$-dimensional case, with steps only in the positive coordinate
directions, was solved by Filaseta~\cite{filaseta}, and a $q$-analogue
was proved by a reflection argument by Krattenthaler~\cite{kratt}.  This
case can be viewed as a variation of the $n$-candidate ballot problem;
how many ways are there to arrange the ballots so that the candidates
stay in order, with the difference between the first-place and
last-place candidates also limited?  (In the ballot problem, equal
coordinates are allowed, but we can translate the start and end by
$(n-1,n-2,\cdots,0)$ to make the inequalities strict.)

Our most important case is the model of $n$ non-colliding particles in
discrete random walks in the interval $(0,m)$.  Equivalently, we can
consider a single random walk in $n$ dimensions in the region
$m>x_1>\cdots>x_n>0$, with each coordinate of the $n$-dimensional walk
corresponding to one of the $n$ particles, and permitted steps in the
positive and negative coordinate directions.  We then simplify the
determinant of infinite sums, computing the exponential generating
function in the number of steps.  This gives the number of walks for one
particle to go from $\eta$ to $\lambda$ while staying in this region, or
for $n$ particles to go from $\eta_i$ to $\lambda_i$ while staying in
the interval $(0,m)$ and not colliding.  The exponential generating
function is
\begin{multline}
\gel(x)=\det_{n\times n}\biggl|
\frac{1}{2m}\sum_{r=0}^{2m-1}
2\sin(\pi r(\lambda_j)/m)
\sin(\pi r(\eta_i)/m)\\
\cdot\exp(2x\cos(\pi r/m))\biggr|.
\end{multline}

Using a reflection argument of Gessel and Krattenthaler~\cite{GK} which
uses the methods of~\cite{KM} and~\cite{HW}, we can also give a formula
for the model of $n$ independent discrete random walks on the circle.
As a one-particle model, this is a variation of the affine Weyl group
$\tilde A_{n-1}$ counting walks to multiple destinations in $\R^n$ which
become equivalent when projected the circle.

\section{Reflectable random walks}

A walk-type is defined by a lattice $L$, a set $S$ of allowable steps
between lattice points, and a region $C$ to which the walks are
confined.  Without affecting the walk problems, we may restrict $L$ and
$C$ to the linear span of $S$, so that $L$, $S$, and $C$ have the same
linear span.  

We will assume $C$ is a {\em Weyl chamber} of a finite Weyl group $W$,
or an {\em alcove} of an affine Weyl group.  The following definitions
and results are given in~\cite{Hu2}.

For a finite Weyl group, we have $L,S,C \subset {\R}^n$; $C$ is
defined by a system of simple roots $\D \subset {\R}^n$ as
\begin{equation}
 C = \{ \vec{x} \in {\R}^n \mid (\alpha, \vec{x})\geq 0 \mbox{ for all }
\alpha \in \D \} ;
\end{equation}
the orthogonal reflections $r_{\alpha}$: $\vec{x} \mapsto \vec{x} - {
\frac{2(\alpha, \vec{x})}{(\alpha,\alpha)}} \alpha$ preserve $L$ and $S$
for all $\alpha$ in $\Delta$; and the $r_{\alpha}$ generate a finite
group $W$ of linear transformations, the {\em Weyl group}.  

The full root system $\Phi$ consists of the images of all roots under
$W$; these roots come in pairs, and we can take the system $\Phi^+$ of
{\em positive roots} to be the set of all roots which are positive
linear combinations of the simple roots; this includes one root from each
pair.  The hyperplanes orthogonal to the positive roots are all
reflections in $W$, and they partition the space into $|W|$ disjoint
Weyl chambers.  

For any root $\alpha$, the corresponding {\em coroot} is
$2\alpha/(\alpha,\alpha)$.  The corresponding affine Weyl group $\tilde
W$ is the semidirect product of $W$ and the translation group of the
coroot lattice; that is, it is generated by reflections
\[
r_{\alpha,k}: \vec{x} \mapsto \vec{x} - {
\frac{2(\alpha, \vec{x})-k}{(\alpha,\alpha)}} \alpha
\]
for all roots $\alpha$ and all integers $k$.  Again, the group $\tilde
W$ contains the reflections $r_{\alpha,k}$ not only for simple $\alpha$
but for all $\alpha$, and these hyperplanes of reflection partition
space into {\em alcoves}.  The alcoves are the regions bounded by
$|\Phi^+|$ simultaneous inequalities
$k_\alpha<(\lambda,\alpha)<k_\alpha+1$, as $\alpha$ runs over all roots
in $\Phi^+$, for $k_\alpha$ in $\Z$.  The {\em principal alcove} $A$ is
bounded by the inequalities $0<(\lambda,\alpha)<1$ for all positive
roots $\alpha$; this can be shown to be non-empty.

{\em Example.} Let $W$ be the symmetric group $S_n$ permuting the $n$
coordinates in $\R^n$; this is generated by reflections in the simple
roots $\D = \{ e_i - e_{i+1}\ ,\ 1\leq i \leq n - 1\}$ (where $e_i$ is
the $i$th coordinate vector), which gives Weyl chamber
$x_1>x_2>\cdots>x_n$.  The positive roots are $e_i-e_j$ for $i<j$; their
hyperplanes $x_i=x_j$ give $n!$ Weyl chambers, one for each permutation
of the coordinates.  The corresponding affine Weyl group contains a
reflection in the hyperplane $x_i-x_j=k$ for any integer $k$; the
principal alcove is thus $x_1>x_2>\cdots>x_n>x_1-1$.  (While this alcove
is unbounded, $W$ really acts only on the subspace in which $\sum
x_i=0$, and the alcove is bounded in that subspace.)

We would like to classify those Weyl chambers and alcoves which allow us
to reflect a walk from the point at which it hits a wall.  The
definition of a reflectable walk from~\cite{decomp} generalizes easily
to the affine case.

{\em Definition.} A walk-type $(L,S,C)$ is {\em reflectable} with
respect to the finite or affine Weyl group $W$ if the steps $S$ are
symmetric under the finite Weyl group, and the following equivalent
conditions 
hold: 
\begin{enumerate}
\item Any step $s \in S$ from any lattice point in the interior of $C$
will not exit $C$. 
\item For each simple root $\alpha_{i}$, there is a real number
$k_i$ such that $(\alpha_i, s) = \pm k_i$ or $0$ for all steps $s\in S$;
$(\alpha_i, \lambda)$ is an integer multiple of $k_i$ for all
$\lambda \in L$; and if $C$ is an alcove of an affine Weyl group, then
$1/k_i$ must be an integer.
\end{enumerate}

To see that these conditions are equivalent, note that if $\alpha_j$ is
any root which is symmetric to $\alpha_i$, then $k_j=k_i$.  Thus the
second condition guarantees than $L$ contains only points whose
dot-product with $\alpha_j$ is an integer multiple of $k_j$, and a
single step can change the dot product 0 or $\pm k_j$.  Therefore, the
walk cannot go from one side of the wall $(\alpha_j,\lambda)=0$ to the
other without stopping on a wall.  Likewise, since each alcove wall is
the set of points with $(\alpha_j,\lambda)=m$ for some $m$, and $m$ is a
multiple of $k_j$, the walk cannot go from one side to the other in a
single step.

{\em Example.} In the example above, in which $W$ is the symmetric
group, the steps $\pm e_i/t$ on $L=\Z^n$ give a reflectable random walk,
with each $k_i=1/t$.  This is reflectable for the affine alcove as well
provided that $t$ is an integer.

By a similar argument, all of the reflectable random walks on the finite
Weyl groups give reflectable walks on the corresponding affine Weyl
groups.  These steps are enumerated in~\cite{decomp}; they turn out to
be precisely the Weyl group images of the {\em minuscule
weights}~\cite{Bou}, those weights with dot-product 0 or $\pm 1$ with
every root.  The reflectable walks include the weights which are
minuscule for only one of $B_n$ and $C_n$.  In the Bourbaki
numbering~\cite{Bou}, the allowed step sets are the Weyl group images of
the following $\Om_i$, the duals of the fundamental roots.

\begin{eqnarray*}
A_n:  &\Om_1,\ldots, \Om_n.   &\mbox{\rm All compatible.} \\
B_n, C_n:  &\Om_1, \Om_n. &\mbox{\rm Not compatible.}\\
D_n:  &\Om_1, \Om_{n-1},\Om_n. &\mbox{\rm All compatible.} \\
E_6:  &\Om_1, \Om_6. &\mbox{\rm Compatible.}\\
E_7:  &\Om_7.\\
E_8,F_4,G_2:  &\mbox{\rm None.}
\end{eqnarray*}

We can also get a reflectable walk by taking any union of compatible
step sets (step sets which give the same $k_i$), and we can add the zero
step.

For the minuscule weights, which include all of the above weights for
the Weyl groups $A_n$, $D_n$, $E_6$ and $E_7$, and also include the
weight $\Om_1$ of $B_n$, and $\Om_n$ of $C_n$, the allowed steps have
dot product 0 or $\pm 1$ with every root; thus multiplying the steps by
$1/t$ for any positive integer $t$ gives a reflectable random walk for
the affine Weyl group.

For the weights $\Om_n$ of $B_n$ and $\Om_1$ of $C_n$, the dot product
with the short roots is 0 or $\pm 1$, and with the long roots is 0 or
$\pm 2$.  Thus, for these cases, we must multiply the steps by $1/2t$
for a positive integer $t$ to satisfy the reflectability condition.

In the natural cases we study later, it will be clear from the
definitions that the first reflectability condition is satisfied.

\section{The reflection argument of Gessel and Zeilberger}

In a reflectable random walk problem, we want to compute $\belk$, the
number of walks from $\eta$ to $\lambda$ of length $k$ which stay in the
interior of a Weyl chamber or alcove.  For example, the ballot problem
can be converted to this form by starting at the point
$(n-1,n-2,\ldots,0)$ instead of the origin, and requiring the
coordinates to remain strictly ordered.  

Let $\cgamk$ denote the number of random walks of length $k$, with steps
in $S$, from the origin to $\gamma$, but {\em unconstrained} by a
chamber.  The fundamental result of Gessel and Zeilberger~\cite{GZ}
(also proved independently by Biane~\cite{biane92} for finite Weyl
groups), is

\begin{thm}\label{gzthm}
If the walk from $\eta$ to $\lambda$ is reflectable, then
\begin{equation}
  \belk=\sumweyl\sgn(w)c_{w(\lambda)-\eta,k}, \label{belkwalk}.
\end{equation}
\end{thm}

If $W$ is an affine group, this is an infinite sum, but only finitely
many terms are nonzero for any fixed $k$.

{\em Proof.} Every walk from $\eta$ to any $w(\lambda)$ which does touch
at least one wall has some first step $j$ at which it touches a wall.
Let the wall be a hyperplane perpendicular to $\alpha_i$, choosing the
largest $i$ if there are several choices~\cite{Pr2}; the reflection in
that wall is a reflection $r_{\alpha_i,k}$.  Reflect all steps of the
walk after step $j$ across that hyperplane; the resulting walk is a walk
from $\eta$ to $r_{\alpha_i,k} w (\lambda)$ which touches the same wall
at step $j$.  This clearly gives a pairing of walks, and since
$r_{\alpha_i,k}$ has sign $-1$, these two walks cancel out
in~\eqref{belkwalk}.  The only walks which do not cancel in these pairs
are the walks which stay within the Weyl chamber or alcove, and since
$w(\lambda)$ is inside the Weyl chamber or alcove only if $w$ is the
identity, this is the desired number of walks. $\qed$

The specific case of the theorem in which $W$ is the symmetric group
$S_n$ was proved by Karlin and McGregor~\cite{KM}.  We can view the
$S_n$ process as $n$ separate walks which are not allowed to collide,
rather than one walk restricted to the region $x_1>\cdots>x_n$, and
interchange two particles when they collide.  If $n_{ij,k}$ is the
number of walks of length $k$ from $\eta_i$ to $\lambda_j$, then the
formula~\eqref{belkwalk} becomes
\begin{equation}
  \belk=\sum_{\sigma\in S_n}\sgn(\sigma) n_{i,\sigma(i),k}
  =\det_{n\times n}\left|n_{ij,k}\right|.
\end{equation}

\section{Unconstrained walks} 

Two natural choices of step sets are the positive and negative
coordinate directions $\pm e_i$, or the $2^n$ diagonals $\pm\1/2 e_1
\cdots\pm\1/2 e_n$.  Both cases give reflectable walks for all the
classical Weyl chambers.  To get a non-trivial walk for the alcoves of
the affine Weyl groups, we have to re-scale either the steps or the
alcoves.  It is more natural to state the problem if we leave the steps
alone and re-scale the alcoves by a factor of $m$ or $2m$, so that we
study the steps $\pm e_i$ in the alcoves of the Weyl group $mW$ or
$2mW$.

For the diagonals, $\cgamk$ is easy to compute, and for the coordinate
directions, the exponential generating function
\[
h_\gamma(x)=\sumkinf\cgamk x^k/k!
\]
is easy to compute.

If the steps $S$ are the diagonals $\pm \1/2 e_1\cdots \pm \1/2 e_n$,
then the walk is essentially $n$ independent walks in the coordinate
directions.  Each step involves a step of $\pm \1/2$ in each coordinate
direction, and thus the walk will go forwards a distance of $\gamma_i$
if there are $\gamma_i$ more forward steps than backward steps.  That
is,
\begin{equation}\label{diagform}
\cgamk = \prodin\binom{k}{(k/2) + \gamma_i}.
\end{equation}

If the steps $S$ are $\pm e_i$, the positive and negative coordinate
directions, let $\chi(\u)=\sum_{i=1}^n {u_i+u_1^{-1}}$, the generating
function for the steps in the formal monomials $u^{(x_1,\ldots,x_n)} =
u_1^{x_1}\cdots u_n^{x_n}$.  Then we have
\[
\cgamk = \chi(\u)^k\Big|_{\u^\gamma},
\]
where $\Big|_{\u^\gamma}$ denotes the coefficient of $\u^\gamma$ in the
polynomial.  This gives 
\[
h_\gamma(x)=\sumkinf\cgamk x^k/k!= \exp(x\chi(\u))\Big|_{\u^\gamma}
\]
for any walk, and for this walk, we have
\[
h_\gamma(x)=\sumkinf\prodin\exp(x(u_i+u_i^{-1}))\Big|_{\u^\gamma}.
\]

This infinite sum can be written as a product of hyperbolic Bessel
functions, using the generating-function definition of the Bessel
functions~\cite{WW}.   We have 
\begin{eqnarray*}
  \exp(x(u+u^{-1})) 
    &=& \sumkinf\frac{x^k}{k!} \sum_{j=-k}^k \binom{k}{j} u^{k-2j}\\
    &=& \sum_{m=-\infty}^{\infty} u^m \sumkinf 
	     \frac{x^k}{k!}\binom{k}{(k+m)/2}\\
    &=& \sum_{m=-\infty}^{\infty} u^m \sumtinf \frac{x^{2t+m}}{t!(t+m)!}\\
    &=& \sum_{m=-\infty}^{\infty} u^m I_m(2x),
\end{eqnarray*}

Thus, in this case, the exponential generating function for the
unconstrained walks is
\begin{equation}\label{coordform}
h_\gamma(x) = \prodin I_{\gamma_i}(2x).
\end{equation}

\section{Non-colliding random walks on an interval, and formulas for
$\tilde C_{\lowercase{n}}$}

The basic techniques are similar for all the classical groups; the case
of $\tilde C_n$ is the simplest as well as the most important case.

The walk with steps $\pm e_i$ in the chamber $m>x_1>x_2>\cdots>x_n>0$
is equivalent to a walk of $n$ independent particles in the interval
$(0,m)$, with the process terminating if two particles collide or if one
particle hits an end of the interval.  It is reflectable if $m$ is an
integer.  The Weyl chamber for $\tilde C_n$ is
$\1/2>x_1>x_2>\cdots>x_n>0$ because one of the roots is $2e_1$.  We thus
rescale the group to $2m\tilde C_n$, to get the chamber
$m>x_1>x_2>\cdots>x_n>0$.

We will also consider the walk for the diagonal steps $\pm \1/2
e_1\cdots \pm \1/2 e_n$.  This is reflectable in the same chamber if $m$
is an integer or half-integer.  

Since the coroots of $\tilde C_n$ are $e_i\pm e_j$ and $\pm e_i$, the
affine Weyl group $2m\tilde C_n$ includes all permutations with any
number of sign changes, and translations of any coordinates by multiples
of $2m$. 

We will first consider the walk for the diagonal steps $\pm \1/2
e_1\cdots \pm \1/2 e_n$.  We write an element of the Weyl group as the
product of $\sigma$ in the symmetric group, reflections
$\epsilon_i=\pm 1$ in the coordinate directions, and translations by
$2mt_i$ in the coordinate directions.
Thus Theorem~\ref{gzthm} gives
\begin{equation}\label{tcnstart}
\belk=\sum_{\sigma\in S_n}\sgn(\sigma)
\sum_{t_i\in\Z}\sum_{\epsilon_i=\pm 1}\prod \epsilon_i
c_{(\epsilon_1\lambda_{\sigma(1)}+2mt_1,\ldots,
    \epsilon_n\lambda_{\sigma(n)}+2mt_n)-\eta,k}.
\end{equation}
Using the formula~\eqref{diagform} gives
\begin{equation}
\belk=\sum_{\sigma\in
S_n}\sgn(\sigma)\sum_{t_i\in\Z}\sum_{\epsilon_i=\pm 1}
\prodin \epsilon_i
\binom{k}{(k/2)+\epsilon_i\lambda_{\sigma(i)}-\eta_i+2mt_i}.
\end{equation}
Splitting each factor into terms with $\epsilon_i=1$ and $\epsilon_i=-1$
gives
\begin{multline}
\belk=\sum_{\sigma\in
S_n}\sgn(\sigma)\\
\sum_{t_i\in\Z}
\prodin \left[\binom{k}{(k/2)+\lambda_{\sigma(i)}-\eta_i+2mt_i}
-\binom{k}{(k/2)-\lambda_{\sigma(i)}-\eta_i+2mt_i}\right].
\end{multline}
We interchange the inner sum and the product to get
\begin{multline}
\belk=\sum_{\sigma\in
S_n}\sgn(\sigma)\\
\prodin \sum_{t_i\in\Z}
\left[\binom{k}{(k/2)+\lambda_{\sigma(i)}-\eta_i+2mt_i}
-\binom{k}{(k/2)-\lambda_{\sigma(i)}-\eta_i+2mt_i}\right].
\end{multline}
And the signed sum over permutations is the definition of a determinant,
which gives
\begin{multline}\label{pretcndiag}
\belk=\\ 
\det_{n\times n}\left|
\sum_{t_i\in\Z} \binom{k}{(k/2)+\lambda_j-\eta_i+2mt_i}
               -\binom{k}{(k/2)-\lambda_j-\eta_i+2mt_i}\right|.
\end{multline}

The infinite periodic sums of binomial coefficients can be written as a
finite sum of powers of roots of unity, and thus of cosines.  Let
$\omega=e^{2\pi i/4m}$; then for fixed $r$ and $s$, we have
\begin{equation}
\omega^{-rs}(\omega^r+\omega^{-r})^k=
\sum_{j}\omega^{r(2j-s)}\binom{k}{(k/2)+j}
\end{equation}
Thus, if we take the sum over all $r$ from 0 to $4m-1$, all terms for
which $2j\not\equiv s \pmod{4m}$ (that is, $j\not\equiv \frac{s}{2}\pmod
{2m}$) will cancel out because of the roots of unity; that is,
\begin{equation}\label{periodicomega}
\frac{1}{4m}\sum_{r=0}^{4m-1}\omega^{-rs}(\omega^r+\omega^{-r})^k=
\sum_{j\equiv \frac{s}{2}\!\!\pmod {2m}} \binom{k}{(k/2)+j}.
\end{equation}
Since the original sum was real, we can eliminate the roots of unity in
this sum by taking the real part of $\omega^{-rs}$ and writing
everything in terms of cosines.  This gives
\begin{equation}\label{periodiccos}
\frac{1}{4m}\sum_{r=0}^{4m-1}\cos(2\pi rs/4m)(2 \cos (2\pi r/4m) )^k=
\sum_{j\equiv \frac{s}{2}\!\!\pmod {2m}} \binom{k}{(k/2)+j}.
\end{equation}

Substituting this formula into~\eqref{pretcndiag} gives a determinant 
formula for this case.
\begin{multline}\label{tcndiag}
\belk=\det_{n\times n}\biggl|
\frac{2^k}{4m}
\sum_{r=0}^{4m-1}
\bigl(\cos(2\pi r 2(\lambda_j-\eta_i)/4m) - 
\cos(2\pi r 2(-\lambda_j-\eta_i)/4m)\bigr)\\
\cdot\cos^k(2\pi r/4m) 
\biggr|.
\end{multline}
We can convert the difference of cosines to a product of sines by the
identity $\cos(\alpha-\beta)-\cos(-\alpha-\beta)=2\sin\alpha\sin\beta$;
doing this and simplifying factors of 2 where possible gives the
simplified formula
\begin{equation}\label{tcndiagsimp}
\belk=\det_{n\times n}\biggl|
\frac{2^{k-1}}{m}
\sum_{r=0}^{4m-1}
\bigl(\sin(\pi r\lambda_j/m)\sin(\pi r \eta_i/m)\bigr)
\cos^k(\pi r/2m) 
\biggr|.
\end{equation}

For $n=1$, this is equivalent to the gambler's ruin
formula~\eqref{onedim} and the similar formula~\eqref{onedimgen}.  The
sum of the stakes $N$ is our $2m$ (and $m$ can be a half-integer), since
our steps are $\pm\1/2$ rather than $\pm 1$.  The gambler's ruin formula
counts the probability a walk will end at 0 on the $k$th step, which is
half the probability it will be at 1 after the $(k-1)$st step, so it
corresponds to walks from $\eta$ to $\lambda=\1/2$ of length $k-1$.  In
addition, we are counting walks rather than determining the
probabilities for $k$-step walks, so we have an extra factor of
$2^{k-1}$.

The final difference is that our sum goes from $r=0$ to $4m-1$ rather
than $2m-1$.  The $r=0$ term is zero.  Adding $2m$ to $r$ changes the
sign of $\cos^k(\pi r/2m)$ if $k$ is odd, and changes the signs of
$\sin(\pi r\lambda/2m)$ (and likewise $\sin(\pi r\eta/2m)$) if $\lambda$
(respectively, $\eta$) is a half-integer.  Thus, if $k$ is odd and
$\lambda-\eta$ is an integer, or $k$ is even and $\lambda-\eta$ is a
half-integer, the terms from $r=2m$ to $4m-1$ cancel out the terms from
$r=0$ to $2m-1$, so the total sum is zero, which is correct because
parity makes such a walk impossible.  Otherwise, the terms are
duplicates, and thus we can multiply by 2 and take only the terms from
$r=0$ to $2m-1$.  The formulas \eqref{onedim} and~\eqref{onedimgen} as
stated are only valid for possible walks, and we must add the extra
terms so that we get zero for impossible walks.

The procedure for the steps $\pm e_i$ is similar.  The same argument
that gave us the determinant~\eqref{pretcndiag} of periodic sums of
binomial coefficients for the case of diagonals shows that the
exponential generating function is a product of periodic sums of Bessel
functions.  We get 
\begin{equation}\label{pretcncoord}
\gel(x)=\sum_{\sigma\in S_n}
\sum_{t_i\in\Z} \sgn(\sigma) \prodin \left(
I_{\lambda_{\sigma(i)}-\eta_i+2mt_i}(2x) -
I_{-\lambda_{\sigma(i)}-\eta_i+2mt_i}(2x)\right). 
\end{equation}
We simplify these periodic sums of Bessel functions by again using the 
generating-function definition for the hyperbolic Bessel functions
\[
\exp(x(z+z^{-1})) = \sum_{j\in \Z}z^j I_j(2x).
\]
If we let $z=\omega^r$, where $\omega=e^{2\pi i/2m}$ here rather than
$e^{2\pi i/4m}$ as for the diagonals, we have
\begin{equation}
\omega^{-rs}\exp(x(\omega^r+\omega^{-r}))=
\sum_{j\in \Z}\omega^{r(j-s)} I_j(2x).
\end{equation}
As before, we now take $1/2m$ times the sum over all values of $r$ to
eliminate the terms for which $j\not\equiv s \pmod{2m}$, and then take
the real part of $\omega^{-rs}$ to write everything in terms of
cosines.  This gives
\begin{equation}\label{besselcos}
\frac{1}{2m}\sum_{r=0}^{2m-1}\cos(2\pi rs/2m)\exp(2x\cos(2\pi r/2m))=
\sum_{j\equiv s\!\! \pmod{2m}} I_j(2x).
\end{equation}

Substituting this into our determinant gives the exponential generating
function for this random walk.
\begin{multline}\label{tcncoord}
\gel(x)=\det_{n\times n}\biggl|
\frac{1}{2m}\sum_{r=0}^{2m-1}
\bigl(\cos(2\pi r(\lambda_j-\eta_i)/2m) -
\cos(2\pi r(-\lambda_j-\eta_i)/2m)\bigr)\\
\cdot\exp(2x\cos(2\pi r/2m))\biggr|.
\end{multline}
As before, we simplify factors of 2 and the difference of cosines to get
the final formula
\begin{multline}\label{tcncoordsimp}
\gel(x)=\det_{n\times n}\biggl|
\frac{1}{m}\sum_{r=0}^{2m-1}
\sin(\pi r \lambda_j/m)
\sin(\pi r \eta_i/m)\\
\cdot\exp(2x\cos(\pi r/m))\biggr|.
\end{multline}

Since this is a model of $n$ non-colliding particles, this formula could
also be derived from the argument of Karlin and McGregor~\cite{KM,KM2}.
The possibility of a walk restricted to an interval is not mentioned
there, but the argument is still valid.  The formula above could thus be
derived from the sum of Bessel functions for the one-dimensional walks
using a particle-interchange argument analogous to the one we use below
in Theorem~\ref{kmhwthm}.

\section{Formulas for $\tilde A_{\lowercase{n}-1}$, and non-colliding
random walks on the circle}

The Weyl group $\tilde A_{n-1}$ gives Weyl chamber
$x_1>x_2>\cdots>x_n>x_1-1$, so we rescale by a factor of $m$ to get
$x_1>x_2>\cdots>x_n>x_1-m$.  The Weyl group acts only on the hyperplane
$\sum x_i=0$.  It is easier to study walks in $\R^n$ with steps in the
coordinate or diagonal directions as before; we will also project these
walks onto the hyperplane on which $\tilde A_{n-1}$ acts to study walks
on this hyperplane.

This Weyl group gives another important reflectable case, the $n$ steps
$e_i$, in addition to the two reflectable cases we had for $\tilde C_n$.
The step set $e_i$ is the primary case that has been studied
previously~\cite{filaseta,kratt}.  This step set does not give a
reflectable walk for $\tilde B_n$, $\tilde C_n$, or $\tilde D_n$ because
the step set is not symmetric under the Weyl group.

In addition, we will study the case of $n$ independent random walks on
the circle (the interval $[0,m]$ with endpoints identified).  We will
use the same notation for this problem as for a single $n$-dimensional
walk, as it is analogous to a single $n$-dimensional walk in $m\tilde
A_{n-1}$.  The desired number of walks $\bpelk$ is the number of ways
for particle $i$ to go from $\eta_i$ to $\lambda_i$ in $k$ steps such
that no two particles collide.  The walk is reflectable if it is
symmetric under an interchange of particles at any time that two
particles are in the same place.

For the steps $e_i$, the problem is the same for the Weyl group $m\tilde
A_{n-1}$ acting on $\R^n$ or the hyperplane, and for $n$ non-colliding
particles on the circle (the interval $[0,m]$ with endpoints identified)
with one moving forward at a time.  The points $\eta$ and
$(\eta_1+c,\ldots,\eta_n+c)$ in $\R^n$ project to the same point on the
hyperplane $\sum x_i=0$, but if there are walks from $\lambda$ to $\eta$
of $k$ steps, then $\sum(\eta_i-\lambda_i)=k$, and there cannot be walks
to $(\eta_1+c,\ldots,\eta_n+c)$ for $c\ne 0$.  (Note that this does not
hold for the other step sets because there are backward steps.)  Thus
there is a bijection between walks in $\R^n$ from $\lambda$ to $\eta$
and walks in the hyperplane from the projection of $\lambda$ to the
projection of $\eta$.  The walks in the hyperplane have steps which are
$(n-1)/n$ in one coordinate direction and $-1/n$ in the other coordinate
directions; that is, they are the Weyl group images of the fundamental
weight $\Om_1$.

To see that the walks on $\R^n$ and of $n$ non-colliding particles on
the circle are equivalent, let each coordinate of a walk in $\R^n$
represent an individual particle.  The first and last particles will
collide if $x_1-m=x_n$, and two adjacent particles will collide if
$x_i=x_{i+1}$; these walls are the same as the walls of the Weyl chamber.
We want to count the number of $k$-step walks in $\R^n$ in which one
particle goes from $\lambda$ to $\eta$.  This is the same as the number
of sets of non-colliding walks of $n$ particles in which particle $i$
goes from $\lambda_i$ to $\eta_i$, provided that the coordinates of
$\eta$ are translated by multiples of $m$ if necessary so that
$\sum(\eta_i-\lambda_i)=k$, and so that they are in decreasing order in
an interval of length $m$.  (Again, this argument is not valid for the
other step sets, because we do not have the condition
$\sum(\eta_i-\lambda_i)=k$.)

For example, a walk of 59 steps from $(2,1,0)$ to $(4,3,5)$ on a circle
of length $m=10$ corresponds to a walk in $\R^n$ ending at $(24,23,15)$;
a walk of 49 steps is impossible without permutation of particles, as
the walk in $\R^n$ cannot end at $(14,23,15)$, and a walk which ends at
$(23,15,14)$ corresponds to a walk on the circle in which the first
particle goes from 2 to 3, not from 2 to 4.  

The number of unconstrained walks from $\eta$ to $\eta+\gamma$ with
steps $e_i$ is the multinomial coefficient
$\frac{k!}{\gamma_1!\cdots\gamma_n!}$.  The Weyl group $m\tilde A_{n-1}$
has coroots $e_i-e_j$, so the Weyl group contains all permutations, with
translations of all coordinates by multiples of $m$ such that the sum of
translations is zero.  Thus, if we look at the walks on $\R^n$ which end
at $\lambda$, then Theorem~\ref{gzthm} gives us the sum for any step set
\begin{equation}\label{tanfree}
\belk=\sum_{\sigma\in S_n}\sum_{\sum t_i=0}
\sgn(\sigma)c_{(mt_1,\ldots,mt_n)+\sigma(\lambda)-\eta,k}.
\end{equation}
which in our case is
\begin{equation}\label{tanfwd}
\belk=\sum_{\sigma\in S_n}\sum_{\sum t_i=0}
\sgn(\sigma)\frac{k!}{\prodin(mt_i+\lambda_{\sigma(i)}-\eta_i)!}
\end{equation}
This is the formula computed by Filaseta and
Krattenthaler~\cite{filaseta,kratt}.

We cannot convert this sum to a single determinant in this case; the
condition $\sum t_i=0$ means that we do not have a periodic infinite
sum.  We can convert the sum to a sum of determinants by interchanging
the order of summation and taking out the constant factor $k!$.  This
gives
\begin{equation}\label{tanfwddet}
\belk=k!\sum_{\sum t_i=0}
\det_{n\times n}\left|1/(mt_i+\lambda_j-\eta_i)!\right|.
\end{equation}

For the other step sets, the problems of $n$ particles on the circle, of
one particle in $\R^n$, and of one particle in the hyperplane $\sum
x_i=0$ on which the Weyl group acts, are not equivalent.  We can use our
methods to get single determinant formulas for walks of $n$ particles in
the circle, or of one particle on the hyperplane.

It is most natural to start with the case of $n$ particles on the
circle, as we will use the results of this case in our formulas for the
walk on the hyperplane.  We will use a reflection argument
from~\cite{KM} closely related to Theorem~\ref{gzthm}, and a technique
used in~\cite{HW} for the analogous problem for Brownian motion.  This
theorem is essentially a case of~\cite[Theorem 2]{GK}.

\begin{thm}\label{kmhwthm}
Let $\cpgamk$ be the number of unconstrained walks in $\R^n$ (not
on the circle) of length $k$ from $\eta$ to $\eta+\gamma$.  We may
assume the coordinates of $\eta$ are in decreasing order; let
$\lambda_s$ be the smallest coordinate of $\lambda$.  If the walk from
$\eta$ to $\lambda$ on the circle of size $m$ is reflectable, then the
number of constrained walks of length $k$ on the circle is
\begin{equation}\label{circlewalk}
\bpelk=\sum_{\sigma\in S_n}\sum_{\sum t_i\equiv s\!\! \pmod{n}}
\sgn(\sigma)c'_{(mt_1,\ldots,mt_n)+\sigma(\lambda)-\eta,k}.
\end{equation}
\end{thm}

{\em Proof.} The walks counted by $\cpgamk$ on $\R^n$ project to walks
of $n$ particles on the circle by taking coordinates modulo $m$.  In a
walk which projects to a good walk, if particle $s$ goes from $\eta_s$
to $\lambda_s+mt_s$, then particle $i$ for $i<s$, which starts at
$\eta_i>\eta_s$, must end between $\lambda_s+mt_s$ and
$\lambda_s+(m+1)t_s$ in order not to collide with particle $s$ when the
walk is projected to the circle.  Since $\lambda_i>\lambda_s$, we must
have $t_i=t_s$.  Similarly, if $i>s$, we have $t_i=t_s-1$.  Thus the sum
of all the $t_i$, corresponding to the total number of revolutions, is
congruent to $s$ modulo $n$.

Now, consider a bad walk counted by~\eqref{circlewalk}.  Consider the
first time at which two particles collide, and the first pair of
particles $i$ and $j$ which collide at that time.  Pair it with the
corresponding walk obtained by switching particles $i$ and $j$ after the
collision.  The new $\sigma$ will differ from the old $\sigma$ by a
transposition; the values of $t_i$ and $t_j$ may change, but the sum
$t_i+t_j$ will not because the total forward distance covered by the
two particles does not change.  Thus the paired walk is still counted
in~\eqref{circlewalk}, and these walks cancel out.

The only walks which do not cancel out are those which have no
collisions, and the only walks which have no collisions and are counted
are those in which the particles end in the correct positions.  $\qed$

We can convert the sum from this theorem to a finite sum in the same way
as before.  For the steps $\pm e_i$, which corresponds to one particle
at a time moving independently, we have the exponential generating
function
\begin{equation}
\gpel(x)=\sum_{\sigma\in S_n}\sum_{\sum t_i\equiv s\!\! \pmod{n}}
\sgn(\sigma)\prodin I_{\lambda_{\sigma(i)}-\eta_i+mt_i}(2x).
\end{equation}
Let $\zeta=e^{2\pi i/mn}$; then we can eliminate the condition ${\sum
t_i\equiv s \pmod{n}}$
by using the fact that
\[
\frac{1}{n}\sum_{u=0}^{n-1}\zeta^{um(-s+\sum t_i)}
\]
is 1 if $\sum t_i\equiv s \pmod{n}$ and 0 otherwise.
Thus we have 
\begin{equation}
\gpel(x)=\frac{1}{n}\sum_{u=0}^{n-1}\zeta^{-ums}\sum_{\sigma\in
S_n}\sum_{t_i \in \Z} \sgn(\sigma)\prodin\zeta^{um t_i}
I_{\lambda_{\sigma(i)}-\eta_i+mt_i}(2x).
\end{equation}
For each value of $u$, we get a determinant as before, but this time
we have period $m$ rather than $2m$.  This gives 
\begin{equation}
\gpel(x)=\frac{1}{n}\sum_{u=0}^{n-1}\zeta^{-ums}\det_{n\times n}\left|
\sum_{t_i\in \Z} \zeta^{um t_i} I_{\lambda_j-\eta_i+mt_i}(2x)\right|.
\end{equation}
Except for the factor $\zeta^{um t_i}$, these determinants are the same
type as as~\eqref{pretcncoord}.  We can again use
$\exp(x(z+1/z))=\sum_{j\in\Z} z^j I_j(2x)$,  and take a sum over
different values of $z$; this time, we will use $z=\zeta^{u+nr}$, and
take the sum of only $m$ terms; that is,  
\begin{equation}
\frac{1}{m}\sum_{r=0}^{m-1}\zeta^{-nrv}\exp(x(\zeta^{u+nr}+\zeta^{-(u+nr)}))=
\sum_{j\equiv v\!\! \pmod{m}} \zeta^{ju} I_j(2x).
\end{equation}

We cannot eliminate all of the complex terms by taking the real part,
but we can get a formula which still contains the complex roots of
unity,
\begin{multline}\label{circlecoord}
\gpel(x)=\frac{1}{n}\sum_{u=0}^{n-1}\zeta^{-ums}\\
\cdot\det_{n\times n}\left|
\frac{1}{m}\sum_{r=0}^{m-1}\zeta^{-nr(\lambda_j-\eta_i)}
\exp(2\cos(2\pi(u+nr)/(mn)))\right|.
\end{multline}

An analogous method works for the diagonal walk.  In this walk on the
circle, the particles all start at integer positions, or all start at
half-integer positions, and at each step, all particles move
simultaneously forwards or backwards by $\1/2$.  We define
$s$ as before and $\zeta=e^{2\pi i/2mn}$ since each periodic sum becomes
$2m$ terms rather than $m$, and follow a similar process to get 
\begin{equation}\label{circlediag}
\bpelk=\frac{1}{n}\sum_{u=0}^{n-1}\zeta^{-2ums}\det_{n\times n}\left|
\frac{2^{k-1}}{m} \sum_{r=0}^{2m-1}
\zeta^{-nr(\lambda_j-\eta_i)}\cos^k(\pi(u+nr)/(mn))
\right|.
\end{equation}

For the case of walks in the $m\tilde A_{n-1}$ Weyl chamber on $\R^n$,
restricted to end at $\eta$, we have the formula~\eqref{tanfree}.  We
have the same difficulty in getting a single determinant for all step
sets, but we can use the same argument as for~\eqref{tanfwddet} to get a
sum of determinants.  For the steps $\pm e_i$, we get the exponential
generating function
\begin{equation}
\gel(x)=\sum_{\sum t_i=0}
\det_{n\times n}\left|I_{mt_i+\lambda_j-\eta_i}(2x)\right|,
\end{equation}
and for the diagonal steps $\pm\1/2 e_1\cdots\pm\1/2e_n$, we get
\begin{equation}
\belk=\sum_{\sum t_i=0}
\det_{n\times n}\left|\binom{k}{(k/2)+mt_i+\lambda_j-\eta_i}\right|.
\end{equation}

In contrast, we can get single determinant formulas for the walk
restricted to the hyperplane $\sum x_i=0$ on which $m\tilde A_{n-1}$
acts, by using the formulas for walks on the circle.  We look at the
walks on $\R^n$, and then project them back to the hyperplane; a walk
which ends at $(\eta_1+c,\ldots,\eta_n+c)$ on $\R^n$ projects to a walk
which ends at $\eta$ on the hyperplane.  In particular, the destinations
$(\eta_1+mt,\ldots,\eta_n+mt)$ for all integers $t$ must all be
considered, and these are exactly the destinations we had in the problem
on the circle.  Summing the formula~\eqref{tanfree} over all such
destinations gives exactly the same sum as in~\eqref{circlewalk},
allowing us to use our previous formulas.  We will get one term for each
$(\eta_1+c+mt,\ldots,\eta_n+c+mt)$ for $0\le c<m$ for which the walk can
reach such destinations.  Thus $c$ may be any integer with $0\le c<m$ if
the steps in $\R^n$ are $\pm e_i$.  It must also be an integer if the
steps are the diagonals, provided that we choose our $\eta$ so that
$\eta_i-\lambda_i-k/2$ is an integer; that is, so that $\eta$ itself has
the correct coordinates to make a walk to $\eta$ possible in $\R^n$ by
parity; if $\eta$ were unreachable by parity, then $c$ would have to be
a half-integer.

For the walk on the hyperplane equivalent to the walk with steps $\pm
e_i$, the $2n$ allowed steps are the projections of $\pm e_i$ on this
hyperplane, which are $(n-1)/n$ in one coordinate direction and $-1/n$
in all others, or $-(n-1)/n$ in one coordinate direction and $1/n$ in
all others.  That is, they are the Weyl group images of the weights
$\Om_1$ and $\Om_{n-1}$.  We define $s$ in~\eqref{circlecoord} as before
(it is not changed when we translate all coordinates of $\eta$ by c),
and take the sum of the $m$ terms to get
\begin{multline}\label{tancoord}
\gel(x)=\sum_{c=0}^{m-1}\frac{1}{n}\sum_{u=0}^{n-1}\zeta^{-ums}\\
\det_{n\times n}\left|
\frac{1}{m}\sum_{r=0}^{m-1}\zeta^{-nr(\lambda_j-\eta_i-c)}
\exp(2\cos(2\pi(u+nr)/(mn)))\right|.
\end{multline}

The walk on $m\tilde A_{n-1}$ corresponding to the diagonals is more
natural.  If a particular diagonal step has $p$ coordinates $\1/2$ and
$n-p$ coordinates $-1/2$, then it projects to the vector in the
hyperplane with $p$ coordinates $(n-p)/n$ and $n-p$ coordinates $-p/n$.
These include all the fundamental weights $\Om_1,\ldots,\Om_{n-1}$ of
the Weyl group $A_{n-1}$, including their images under permutations; we
also get the zero step with multiplicity 2 from the two diagonals with
all coordinates equal.  
\begin{multline}\label{tandiag}
\belk=\sum_{c=0}^{m-1}\frac{1}{n}\sum_{u=0}^{n-1}\zeta^{-2ums}\\
\cdot\det_{n\times
n}\left| 
\frac{2^{k-1}}{m} \sum_{r=0}^{2m-1}
\zeta^{-nr(\lambda_j-\eta_i-c)}\cos^k(\pi(u+nr)/(mn))
\right|.
\end{multline}

\section{Formulas for $\tilde B_n$ and $\tilde D_n$}

The affine Weyl group $2m\tilde B_n$ gives the chamber
$x_1>x_2>\cdots>x_n>0, x_1+x_2<2m$.  The affine Weyl group $2m\tilde
D_n$ gives the chamber $x_1>x_2>\cdots>x_n, x_1+x_2<2m, x_{n-1}>-x_n$.
Neither one gives a natural model for $n$ independent particles in one
dimension.  

The Weyl group $2m\tilde B_n$ contains all permutations with any number
of sign changes, and with translations of coordinates by multiples of
$2m$ such that the total translation is a multiple of $4m$ (since
$\tilde B_n$ does not have $2e_i$ as a root and thus does not have
$e_i$ as a coroot).  The Weyl group $2m\tilde D_n$ contains all
permutations with an even number of sign changes, and with translations
of coordinates by multiples of $2m$ such that the total translation is a
multiple of $4m$.

We will find the formulas for these cases by using the fact that
$2m\tilde B_n$ and $2m\tilde D_n$ are subgroups of $2m\tilde C_n$, of
index 2 and 4.  We can thus modify the formula~\eqref{tcnstart} by
including a term which is 1 if the potential endpoint of the walk is in
the Weyl group image of the appropriate subgroup, and 0 otherwise; this
is analogous to the technique used to get $D_n$ formulas from $B_n$
formulas in~\cite{brownian,decomp}.

If an element of $2m\tilde C_n$ contains translations by $2mt_i$, then it
is in $2m\tilde B_n$ if $\sum t_i$ is even, and thus
$\bigl(1+(-1)^{\sum t_i}\bigr)/2$ is a factor which is 1 if the element is
in $2m\tilde B_n$ and 0 otherwise.  Likewise, if it contains reflections
$\epsilon_i$, it is only in $2m\tilde D_n$ if it is in $2m\tilde B_n$
and $\prod \epsilon_i=1$, so the appropriate factor is $\bigl(1+\prod
\epsilon_i\bigr)/2$.

Thus, for $2m\tilde B_n$, we have
\begin{multline}\label{tbnstart}
\gel(x)=\\
\sum_{\sigma\in S_n}\sum_{\epsilon_i=\pm 1}
\sum_{t_i\in\Z} \frac{1+(-1)^{\sum t_i}}{2} \sgn(\sigma) 
\prodin \epsilon_i I_{\epsilon_i\lambda_{\sigma(i)}-\eta_i+2mt_i}(2x). 
\end{multline}
We take the $1/2$ and $(-1)^{\sum t_i}/2$ terms separately.  The $1/2$
term is the term for $2m\tilde C_n$.  The other term also gives a
determinant, 
\begin{equation}
\det_{n\times n}\left| \sum_{t_i\in\Z} (-1)^{t_i} \left(
I_{\lambda_j-\eta_i+2mt_i}(2x) -
I_{-\lambda_j-\eta_i+2mt_i}(2x)\right)\right|.
\end{equation}
These sums are periodic, but with period $4m$ rather than $2m$.  We can
express the terms with odd and even $t_i$ as separate periodic sums as
before,
\begin{multline}
\frac{1}{4m}\sum_{u=0}^{4m-1}
\bigl(\cos(2\pi us/4m)-\cos(2\pi u(2m+s)/4m)\bigr)
\exp(2x\cos(2\pi u/4m))\\
=\sum_{j\equiv s\!\! \pmod{4m}} I_j(2x) - I_{2m+j}(2x).
\end{multline}
For even $u$, $\cos(2\pi us/4m)=\cos(2\pi u(2m+s)/4m)$, so these terms
are all zero.  For odd $u$, $\cos(2\pi us/4m)=-\cos(2\pi u(2m+s)/4m)$, so
the two cosines combine to one term.  Thus we can let $u=2r+1$ and write
the sum as 
\begin{multline}
\frac{1}{2m}\sum_{r=0}^{2m-1}
\cos(2\pi (2r+1)s/4m)
\exp(2x\cos(2\pi (2r+1)/4m))\\
= \sum_{j\equiv s\!\! \pmod{4m}} I_j(2x) - I_{2m+j}(2x).
\end{multline}
And we put both of these terms together in the determinant, to get
\begin{multline}\label{tbncoord}
\gel(x) = \frac{1}{2}\Biggl[\det_{n\times n}\biggl|
\frac{1}{2m}\sum_{r=0}^{2m-1}
\Bigl[\cos\Bigl(\frac{2\pi r(\lambda_j-\eta_i)}{2m}\Bigr) -
\cos\Bigl(\frac{2\pi r(-\lambda_j-\eta_i)}{2m}\Bigr)\Bigr]\\
\cdot\exp(2x\cos(2\pi r/2m))\biggr|\\
 +\det_{n\times n}\biggl|
\frac{1}{2m}\sum_{r=0}^{2m-1}
\Bigl[\cos\Bigl(\frac{2\pi (2r+1)(\lambda_j-\eta_i)}{4m}\Bigr) -
\cos\Bigl(\frac{2\pi (2r+1)(-\lambda_j-\eta_i)}{4m}\Bigr)\Bigr]\\
\cdot \exp(2x\cos(2\pi (2r+1)/4m))\biggr|\Biggr].
\end{multline}
Again, we simplify this to
\begin{multline}\label{tbncoordsimp}
\gel(x) = \frac{1}{2}\Biggl[\det_{n\times n}\biggl|
\frac{1}{m}\sum_{r=0}^{2m-1}
\sin(\pi r \lambda_j/m)
\sin(\pi r \eta_i/m)\\
\cdot\exp(2x\cos(\pi r/m))\biggr|\\
 +\det_{n\times n}\biggl|
\frac{1}{m}\sum_{r=0}^{2m-1}
\sin(\pi (2r+1) \lambda_j/2m)
\sin(\pi (2r+1) \eta_i/2m)\\
\cdot \exp(2x\cos(\pi (2r+1)/2m))\biggr|\Biggr].
\end{multline}

Likewise, for $2m\tilde D_n$, we have as our initial formula
\begin{multline}\label{tdnstart}
\gel(x)=\sum_{\sigma\in S_n}\sum_{\epsilon_i=\pm 1}
\sum_{t_i\in\Z} \frac{1+(-1)^{\sum t_i}}{2} \frac{1+\prod \epsilon_i}{2}
\sgn(\sigma)\\
\prodin I_{\epsilon_i\lambda_{\sigma(i)}-\eta_i+mt_i}(2x). 
\end{multline}
We split this into four terms, taking each combination of the $\1/2$ or
the other term.  When we use $\1/2$ rather than $\prod \epsilon_i/2$,
and continue as in~\eqref{tcnstart} or~\eqref{tbnstart}, we have a plus
sign rather than a minus sign before the term of
$I_{-\lambda_{\sigma(i)}-\eta_i+mt_i}(2x)$.  Everything else carries out
just as before; we have a plus sign rather than a minus sign between the
two terms in this determinant.  This
gives us the determinant formula for 
$2m\tilde D_n$,
\nopagebreak[3] 
\begin{multline}\label{tdncoord}
\gel(x)=\\
\frac{1}{4}\Biggl[\det_{n\times n}\biggl|
\frac{1}{2m}\sum_{r=0}^{2m-1}
\Bigl[\cos\Bigl(\frac{2\pi r(\lambda_j-\eta_i)}{2m}\Bigr) -
\cos\Bigl(\frac{2\pi r(-\lambda_j-\eta_i)}{2m}\Bigr)\Bigr]\\
\shoveright{\cdot\exp(2x\cos(2\pi r/2m))\biggr|}\\
+\det_{n\times n}\biggl|
\frac{1}{2m}\sum_{r=0}^{2m-1}
\Bigl[\cos\Bigl(\frac{2\pi (2r+1)(\lambda_j-\eta_i)}{4m}\Bigr) -
\cos\Bigl(\frac{2\pi (2r+1)(-\lambda_j-\eta_i)}{4m}\Bigr)\Bigr]\\
\shoveright{\cdot\exp(2x\cos(2\pi (2r+1)/4m))\biggr|}\\
+\det_{n\times n}\biggl|
\frac{1}{2m}\sum_{r=0}^{2m-1}
\Bigl[\cos\Bigl(\frac{2\pi r(\lambda_j-\eta_i)}{2m}\Bigr) +
\cos\Bigl(\frac{2\pi r(-\lambda_j-\eta_i)}{2m}\Bigr)\Bigr]\\
\shoveright{\cdot\exp(2x\cos(2\pi r/2m))\biggr|}\\
+\det_{n\times n}\biggl|
\frac{1}{m}\sum_{r=0}^{2m-1}
\Bigl[\cos\Bigl(\frac{2\pi (2r+1)(\lambda_j-\eta_i)}{4m}\Bigr) +
\cos\Bigl(\frac{2\pi (2r+1)(-\lambda_j-\eta_i)}{4m}\Bigr)\Bigr]\\
\exp(2x\cos(2\pi (2r+1)/4m))\biggr| \Biggr].
\end{multline}
We simplify as before, with the sum of cosines simplifying by 
$\cos(\alpha-\beta)+\cos(-\alpha-\beta)=2\cos\alpha\cos\beta$, to get our
final formula for $2m\tilde D_n$
\nopagebreak[3] 
\begin{multline}\label{tdncoordsimp}
\gel(x)=\\
\frac{1}{4}\Biggl[\det_{n\times n}\biggl|
\frac{1}{m}\sum_{r=0}^{2m-1}
\sin(\pi r \lambda_j/m)
\sin(\pi r \eta_i/m)\\
\shoveright{\cdot\exp(2x\cos(\pi r/m))\biggr|}\\
+\det_{n\times n}\biggl|
\frac{1}{m}\sum_{r=0}^{2m-1}
\sin(\pi (2r+1) \lambda_j/2m)
\sin(\pi (2r+1) \eta_i/2m)\\
\shoveright{\cdot\exp(2x\cos(\pi (2r+1)/2m))\biggr|}\\
+\det_{n\times n}\biggl|
\frac{1}{m}\sum_{r=0}^{2m-1}
\cos(\pi r \lambda_j/m)
\cos(\pi r \eta_i/m)\\
\shoveright{\cdot\exp(2x\cos(\pi r/m))\biggr|}\\
+\det_{n\times n}\biggl|
\frac{1}{m}\sum_{r=0}^{2m-1}
\cos(\pi (2r+1) \lambda_j/2m)
\cos(\pi (2r+1) \eta_i/2m)\\
\exp(2x\cos(\pi (2r+1)/2m))\biggr| \Biggr].
\end{multline}

In both cases, we have analogous formulas for the diagonal walk, which
again allows $m$ to be a half-integer.

\section{Open questions}

The number of reflectable random walks in a classical Weyl chamber is of
interest in representation theory.  Let the step set be the set of
weights of a representation of the corresponding Lie group.  Let the
starting point $\eta$ be $\rho$, half the sum of the positive roots, and
the end point $\lambda$ be $\rho+\mu$.  If the walk is reflectable, then
the number of walks from $\rho$ to $\rho+\mu$ is the multiplicity of the
representation with highest weight $\mu$ in the $k$th tensor power of
the representation whose weights are the step set~\cite{decomp}.  Does
the number of walks in an alcove of an affine Weyl chamber have a
similar meaning in representation theory, either in the representations
of other Lie groups or Lie algebras, or of the classical Lie groups over
some other field?

We have a formula for the probability that $n$ particles in an interval
or on a circle will not collide in $k$ steps.  Can we get a single
general formula, or an asymptotic, for the total probability that there
will be no collision?

There are formulas and interpretations of the $q$-analogue of the case
with Weyl group $\tilde A_{n-1}$ and steps $+e_i$~\cite{kratt}.  Can
$q$-analogues of the other cases be defined, and are they of interest?

\end{document}